\documentclass[12pt]{article}

\usepackage{amsmath}
\usepackage{amssymb}

\def\la{\lambda}

\newtheorem{theorem}{Theorem}

\newcommand{\StirlingPartition}[2]{\genfrac{ \{ }{ \} }{0pt}{}{#1}{#2}}
\newcommand{\StirlingCycle}[2]{\genfrac{[}{]}{0pt}{}{#1}{#2}}

\newcommand{\un}[1]{\ensuremath{\underline{#1}}}

\begin{document}
\begin{center}
{\Large Jordan and Smith forms of Pascal-related matrices \\ 
}
\vspace{10mm}
DAVID CALLAN  \\
Department of Statistics  \\
University of Wisconsin-Madison  \\
1210 W. Dayton St   \\
Madison, WI \ 53706-1693  \\
{\bf callan@stat.wisc.edu}  \\
\vspace{5mm}
\today
\end{center}

\vspace{5mm}

%\vspace*{5mm}

%\emph{Keywords:} Pascal matrix; Jordan form; Smith normal form; Stirling numbers

\vspace*{5mm}

\textbf{1.   Definitions}  \quad
Let $P_{n}=\left(\binom{i-1}{j-1}\right)_{1 
\le i,j \le n}$, $S_{n}=\left(\StirlingPartition{i-1}{j-1}\right)_{1
\le i,j \le n}$, $C_{n}=\left(\StirlingCycle{i-1}{j-1}\right)_{1
\le i,j \le n}$ denote the $n$-by-$n$ lower triangular 
matrices of binomial coefficients, Stirling partition numbers and 
Stirling cycle numbers respectively. For example,

\[
P_{5}=
\left(\begin{array}{ccccc}
	1 &  &  &  &   \\
	1 & 1 &  &  &   \\
	1 & 2 & 1 &  &   \\
	1 & 3 & 3 & 1 &   \\
	1 & 4 & 6 & 4 & 1
\end{array}\right),\ \ 
S_{5}=
\left(\begin{array}{ccccc}
	1 &  &  &  &   \\
	0 & 1 &  &  &   \\
	0 & 1 & 1 &  &   \\
	0 & 1 & 3 & 1 &   \\
	0 & 1 & 7 & 6 & 1
\end{array}\right),\ \  
C_{5}=
\left(\begin{array}{ccccc}
	1 &  &  &  &   \\
	0 & 1 &  &  &   \\
	0 & 1 & 1 &  &   \\
	0 & 2 & 3 & 1 &   \\
	0 & 6 & 11 & 6 & 1
\end{array}\right).
\]

We use $i^{\underline{j}}$ and  $i^{\overline{j}}$  for the
falling and rising factorials respectively and we 
adopt the convention a binomial coefficient is 
zero if either of its parameters is negative.
Let
$F_{n,r}$ denote the $n$-by-$n$ matrix 
$\Big((i-1)^{\underline{i-j}}\binom{i-j-1}{r-1}\Big)_{1 \le 
i,j \le n}$. 
Note that $F_{n,r}$ is a block 
matrix $\left(\begin{smallmatrix} 0 & 0\\
G_{n,r} & 0 \end{smallmatrix}\right)$ with $G_{n,r}$ 
lower triangular of size $n-r$.
Let $H_{n,r}$ denote the $(n-r)$-by-$(n-r)$ lower triangular banded 
matrix
$\left( 
(-1)^{(i-j)}\binom{i-1}{j-1} r^{\underline{i-j}}\right)_{1 \le i,j \le 
n-r}$. Note the $(i,j)$ entry of $H_{n,r}$ is 0 if $i-j>r$  and 
can also be expressed as $
\binom{i-1}{j-1}(-r)^{\overline{i-j}}$.
 Let $D_{n,r}$ denote the diagonal 
matrix with the same diagonal, $(i^{\overline{r}})_{1 \le i \le 
n-r}$, as $G_{n,r}$.
 For example,
\[
G_{6,2}=
\left(\begin{array}{cccc}
	2 &  &  &     \\
	12 & 6 &  &     \\
	72 & 48 & 12 &     \\
	480 & 360 & 120 & 20  
\end{array}\right),\ \ 
H_{6,2}=
\left(\begin{array}{cccc}
	1 &  &  &     \\
	-2 & 1 &  &     \\
	2 & -4 & 1 &     \\
	0 & 6 & -6 & 1    
\end{array}\right),\ \  
D_{6,2}=
\left(\begin{array}{cccc}
	2 &   &     &     \\
      & 6 &     &     \\
	  &   & 12  &     \\
	  &   &     & 20   
\end{array}\right).
\] 

\vspace*{5mm}

\textbf{2.  Matrix Identities}  \quad
With the preceding definitions, we have the following four identities:
\begin{equation}
	S_{n}^{-1}\,P_{n}\,S_{n}=\left(
	\begin{array}{ccccc}
		1 &  &  &  &   \\
		1 & 1 &  &  &   \\
		 & 2 & 1 &  &   \\
		 &  & \ddots & \ddots &   \\
		 &  &  & n-1 & n
	\end{array}\right)
	\label{eq:1}
\end{equation}
\begin{equation}
	\left(\StirlingCycle{i}{j}\right)_{1\le i,j \le n}
	\left(\binom{i}{j}\right)_{1\le i,j \le n}
	\left(\StirlingCycle{i}{j}\right)_{1\le i,j \le n}^{-1}=
	\text{diag}\Big((1,2,\ldots,n)\Big)
\end{equation}
\begin{equation}
	C_{n}\,(P_{n}-I_{n})^{r}\,C_{n}^{-1}=F_{n,r}
	\label{eq:3}
\end{equation}
\begin{equation}
	G_{n,r}\,H_{n,r}=D_{n,r}
\end{equation}

\vspace*{5mm}

 \textbf{3. Canonical Forms}   \quad
Identity (1) puts $P_{n}$ in bidiagonal form, which is good enough to 
obtain the Jordan form: say a matrix is  
near-Jordan if it has the form
\[
J_{n}(\la;c_{1},\ldots,c_{n-1}) =
\left(
	\begin{array}{ccccc}
		\la &  &  &  &   \\
		c_{1} & \la &  &  &   \\
		 & c_{2} & \la &  &   \\
		 &  & \ddots & \ddots &   \\
		 &  &  & c_{n-1} & \la
	\end{array}\right)
\]
with all the $c_{i}$'s nonzero. (A Jordan matrix is one with each 
$c_{i}=1$.) Now $J_{n}(\la;c_{1},\ldots,c_{n-1})$ is similar to a Jordan 
matrix via 
\[
D\, J_{n}(\la;c_{1},\ldots,c_{n-1})\, D^{-1} = J_{n}(\la;1,\ldots,1)
\]
where $D= \text{diag}\Big( (1,c_{1},c_{1}c_{2},\ldots,c_{1}c_{2}\cdots 
c_{n-1})\Big)$. Hence the Jordan form of $P_{n}$ consists of a single 
block with $\la=1$. We can deduce the Jordan form of $P_{n}$ mod $p$: 
reducing (1) mod $p$, the right side becomes block diagonal with as 
many $p$-by-$p$ near-Jordan blocks as will fit in an $n$-by-$n$ 
matrix and (possibly) one of smaller size. Hence we have 
\begin{theorem}
	The Jordan form of 
    $P_{n}$ mod $p$ has precisely $\lceil \frac{n}{p} \rceil$ blocks $($all 
    corresponding to the eigenvalue 1$)$ 
    and $\lfloor \frac{n}{p} \rfloor$ of them are of size $p$. 
\end{theorem}
In particular, for $n\ge p$  the minimal polynomial of $P_{n}$ 
mod $p$ is $(x-1)^{p}$ \cite{bayat}.

Identity (2) shows
\begin{theorem}
	The matrix $\Big(\StirlingCycle{i}{j}\Big)_{1 \le i,j \le n}^{-1} = 
    \Big((-1)^{i-j}\StirlingPartition{i}{j}\Big)_{1 \le i,j \le n}$ 
    is a matrix of eigenvectors for 
    $\Big( \binom{i}{j}\Big)_{1\le i,j \le n}$. 
\end{theorem}

Recall that two integer matrices $A,B$ are \emph{equivalent} if there exist unimodular 
matrices $U,V$ such that $UAV=B$.
From (4), we have 
$\left(
\begin{smallmatrix} 0 & I_{n-r} \\
I_{r} & 0
\end{smallmatrix}\right)
\left(
\begin{smallmatrix} 0 & 0 \\
G_{n,r} & 0
\end{smallmatrix}\right)
\left(
\begin{smallmatrix} H_{n,r} & 0 \\
0 & 0
\end{smallmatrix}\right) =
\left(
\begin{smallmatrix} D_{n,r} & 0 \\
0 & 0
\end{smallmatrix}\right)$ where the middle matrix in the product is 
$F_{n,r}$. Combined with (3), we have
\begin{theorem}
	$(P_{n}-I_{n})^{r}$ is 
    equivalent to the diagonal matrix 
    diag$(1^{\overline{r}},2^{\overline{r}},\ldots,(n-r)^{\overline{r}},0,\ldots,0)$.
\end{theorem}
This is good enough \cite[Theorem II.13, p.\,30]{integral matrices} to obtain the 
elementary divisors and hence the Smith normal 
form  of
$(P_{n}-I_{n})^{r}$.

\vspace*{5mm}

 \textbf{4. Proofs}   \quad
After left multiplying by $S_{n}$, equating 
the  $(n+1,m+1)$ entries and using the defining recurrence 
$\StirlingPartition{n+1}{m+1}=
\StirlingPartition{n}{m}+(m+1)\StirlingPartition{n}{m+1}$, 
identity (1) reduces to
$\StirlingPartition{n+1}{m+1}=\sum_{k}\binom{n}{k}\StirlingPartition{k}{m}$
\cite[6.15,\ p.\,265]{gkp}. 
Similarly, identity (2) is equivalent to a simple variant of
$
\StirlingCycle{n+1}{m+1}=\sum_{k}\StirlingCycle{n}{k}\binom{k}{m}
$
\cite[6.16,\ p.\,265]{gkp} 

Equating the  $(n,m)$ entries, (3) is equivalent to
\[
\sum_{k=m+r}^{n}\StirlingCycle{n}{k}\binom{k}{m}\ 
r!\StirlingPartition{k-m}{r}= \sum_{k=m}^{n-r}n^{\underline{n-k}} 
\binom{n-k-1}{r-1}\StirlingCycle{k}{m}.
\]
Both sides count the following set $S$ of combinatorial objects: an 
object of $S$ is a partition of $[n]$ into cycles, $m$ of which are 
colored black (say) and the rest are colored from $\{1,2,\ldots,r\}$ 
so that all $r$ colors appear. The left side counts $S$ by number of 
cycles: $\StirlingCycle{n}{k}$---partition $[n]$ into $k$ cycles; 
$\binom{k}{m}$---choose the black cycles; 
$r!\StirlingPartition{k-m}{r}$---assign colors to the other $k-m$ 
cycles. The right side counts $S$ by total number of elements in the 
black cycles: $n^{\un{n-k}}$\:---form a list of $n-k$ elements from 
$[n]$; $\binom{n-k-1}{r-1}$---choose $r-1$ of the $n-k-1$ ``dividers'' 
separating elements of the list, thereby obtaining a list of $r$ 
nonempty lists each of which is canonically converted to cycles 
(for example, canonical might mean that cycles start 
at the  left to right minima  of a list), all 
cycles from the $i$th sublist getting color $i$; 
$\StirlingCycle{k}{m}$---choose $m$ black cycles on the remaining $k$ elements.

Identity (4) is a consequence of the binomial theorem for rising 
factorials which implies that the inverse of $H_{n,r}=
\left( 
\binom{i-1}{j-1} (-r)^{\overline{i-j}}\right)_{1 \le i,j \le 
n-r}$ is $\left( 
\binom{i-1}{j-1} r^{\overline{i-j}}\right)_{1 \le i,j \le 
n-r}$.

\vspace*{5mm}

\textbf{5. Further Remarks}   \quad
The $(i,j)$ entry of $(P_{n}-I_{n})^{r}$ has a closed form, $r! 
\StirlingPartition{i-j}{r}\binom{i-1}{j-1}$, 
for which there
is a simple combinatorial explanation. For a sequence 
$\mathbf{c}=(c_{i})_{i\ge 0}$, set 
$P_{n}(\mathbf{c})=\left(c_{i-j}\binom{i-1}{j-1}\right)_{1\le i,j \le n}$. 
Then $P_{n}(\mathbf{c})P_{n}(\mathbf{d})=P_{n}(\mathbf{c}*\mathbf{d})$ 
where 
$\mathbf{c}*\mathbf{d}=\left(\sum_{i=0}^{n}\binom{n}{i}c_{i}d_{n-i}\right)_{n\ge 0}$.
If $\mathbf{c}$ counts labeled structures of some kind ($c_{n}$ is 
the number of structures on $[n]$), then the $r$-fold convolution 
$\mathbf{c}*\ldots*\mathbf{c}$ counts $r$-sequences of structures 
whose labels partition $[n]$. In particular, for (the species of) 
nonempty sets  $\mathbf{c}=(0,1,1,1,\ldots)$ and 
$\mathbf{c}*\ldots*\mathbf{c}$ counts ordered partitions of $[n]$ 
into $r$ (unordered  nonempty) blocks; the number of such ordered 
partitions is 
$r!\StirlingPartition{n}{r}$. Hence, with 
$\mathbf{c}:=(0,1,1,1,\ldots)$ and 
$\mathbf{d}:=\left(r!\StirlingPartition{i}{r}\right)_{i\ge 0},\ 
(P_{n}-I_{n})^{r}=P_{n}(\mathbf{c})^{r}=
P_{n}(\mathbf{c}*\ldots * \mathbf{c})=P_{n}(\mathbf{d})=
\Big(r! 
\StirlingPartition{i-j}{r}\binom{i-1}{j-1}\Big)_{1 \le i,j \le n}$, as 
asserted.

Say an integer matrix is equivalent to its diagonal if it is equivalent 
to the matrix obtained by zeroing out all offdiagonal entries. Suppose 
a sequence $\mathbf{c}$ (such as $\left(\StirlingPartition{i}{r}\right)_{i\ge 0}$) 
begins with $r$ 0's so that $P_{n}(\mathbf{c})$ has the block form $
\left(
\begin{smallmatrix} 0 & 0 \\
Q_{n}(\mathbf{c}) & 0
\end{smallmatrix}\right)$ with $Q_{n}(\mathbf{c})$ square of size 
$n-r$. It follows from identities (3) and (4) 
above that the following statement is 
true. For $\mathbf{c}=\left(\StirlingPartition{i}{r}\right)_{i\ge 
0},\ Q_{n}(\mathbf{c})$ is equivalent to its diagonal. It seems the 
same statement is true for $\mathbf{c}=\left(\StirlingCycle{i}{r}\right)_{i\ge 
0}$ and it would be interesting to have a proof.

\end{document}